\documentclass[12pt]{amsart}
\usepackage{amscd,enumerate,amsfonts,calc,amsmath,verbatim}

\begin{document}
\title{Chains and families of tightly closed ideals}
\author{Adela Vraciu}

\date{\today}
\thanks{ }
\begin{abstract} 
We prove tight closure analogues of results of Watanabe about chains and families of integrally closed ideals.
\end{abstract}

\maketitle

\swapnumbers
\theoremstyle{plain}
\newtheorem{theorem}{Theorem}[section]

\newtheorem{prop}[theorem]{Proposition}
\newtheorem{lemma}[theorem]{Lemma}
\newtheorem{corollar}[theorem]{Corollary}
\newtheorem*{Corollary}{Corollary}

\theoremstyle{definition}
\newtheorem{note}[theorem]{Note}
\newtheorem{obs}[theorem]{Observation}
\newtheorem{definition}[theorem]{Definition}
\newtheorem*{Definition}{Definition}
\newtheorem{example}[theorem]{Example}
\newtheorem*{notation}{Notation}
\newtheorem*{conj}{Conjecture}
\newtheorem*{claim}{Claim}
\newtheorem*{question}{Question}

\newcommand{\li}{\tilde}
\newcommand{\aaa}{\mathfrak{a}}
\newcommand{\bbb}{\mathfrak{b}}
\newcommand{\ccc}{\mathfrak{c}}
\newcommand{\ub}{\underline{b}}

\newcommand{\m}{\mathfrak{m}}
\newcommand{\param}{\underline{x}}
\newcommand{\tpar}{\underline{x}^{[t]}}
\newcommand{\tparq}{\underline{x}^{[tq]}}
\newcommand{\bs}{\boldsymbol}
\newcommand{\tx}{\noindent \textbf}
\newcommand{\ld}{\ldots}
\newcommand{\cd}{\cdots}
\newcommand{\q}{^{[q]}}
\newcommand{\cor}{^{<q>}}
\newcommand{\spec}{I^{*sp}}
\newcommand{\f}{(f_1, \ld, \hat{f_i}, \ld, f_n)}
\newcommand{\eqq}{\Leftrightarrow}
\newcommand{\fs}{(f', f_1, \ld, \hat{f_i}, \ld , f_n)}
\newcommand{\ins}{I_1\cap \ld \cap \hat{I_i} \cap \ld \cap I_n}
\newcommand{\xs}{x_1, \ld, x_d}
\newcommand{\xt}{(x_1^t, \ld, x_d^t)}
\newcommand{\xtq}{(x_1^{tq}, \ld, x_d^{tq})}
\newcommand{\product}{x_1\cdots x_d}

\newcommand{\ui}{\underline{i}}
\newcommand{\und}{\underline}
\newcommand{\arr}{\Rightarrow}
\newcommand{\lar}{\longrightarrow}
\newcommand{\inc}{\subseteq}
\newcommand{\w}{\omega}

\newcommand{\Tor}{\mathrm{Tor}^R}
\newcommand{\Ext}{\mathrm{Ext}_R}
\newcommand{\x}{\underline{x}}
\newcommand{\ol}{\overline}
\newcommand{\syz}{\mathrm{syz}}
\newcommand{\Soc}{\mathrm{Soc}}

\section*{Introduction}
In \cite{W}, Watanabe proved the following result:
\begin{theorem}(\cite{W}, Theorem 2.1)
Let $(A, \m)$ be an excellent normal local ring with algebraically closed residue field, and let $J\subset I$ be $\m$-primary integrally closed ideals. 
Then there exist integrally closed ideals $I'$ with $J\subset I' \subset I$ and $\lambda (I/I')=1$.
 \end{theorem}
(Here, and in the rest of this paper, $\lambda$ denotes length.)

As an immediate corollary, one sees that there exists a sequence $J=I_0\subset I_1 \ldots \subset I_n=I$ consisting of integrally closed ideals, such that $\lambda(I_{i+1}/I_i)=1$ for all $i$ (\cite{W}, Corollary 2.2).

One of the main results of this paper (Corollary ~\ref{chain}) is that the corresponding statement holds for tightly closed ideals.

 Note that we are assuming that the residue field is perfect instead of algebraically closed, and we do not need the assumption 
that $I$ and $J$ are $\m$-primary. 

Watanabe also shows that the family of adjacent integrally closed ideals ``looking down from above'', i.e. the family of 
integrally closed ideals $I'\subset I$ with $\lambda(I/I')=1$ (where $I$ is a fixed $\m$-primary integrally closed ideal) 
is in one-to-one correspondence with points of an algebraic variety of dimension $d-1$, where $d$ is the Krull dimension of the ring. Since $I$ is $\m $-primary, $d$ is also the analytic spread of $I$.

The corresponding result for tightly closed ideals (Corollary~\ref{gras}) is proved in Section 2 of this paper. It shows that the 
family of tightly closed adjacent ideals ``looking down from above'', i.e. the tightly closed ideals $I'\subset I$ 
with $\lambda(I/I')=1$ (where $I$ is a fixed tightly closed ideal) are in one-to-one correspondence with points of an $l-1$ dimensional vector space, where $l$ is the $*$-spread of $I$ (see Observation ~\ref{spread} for the definition of $*$-spread). Since the $*$-spread of an ideal is the tight closure analog of the analytic spread, Corollary ~\ref{gras} is an exact analog of Watanabe's result.

The methods employed here allow one to recover analogous results for any closure operation satisfying certain axiomatic conditions, namely: 
a version of Nakayama's lemma (Proposition ~\ref{nak} is the tight closure notion); a notion of ``spread''  (if $I^{cl}$ denotes the closure of $I$, we say that $J$ is a minimal $cl$-reduction of $I$ if $I\subset J^{cl}$, and $J$ is minimal with this property with respect to inclusion; in order for the closure $cl$ to admit spread, all minimal reductions of a given ideal $I$ must have the same minimal number of generators, called the $cl$-spread of $I$); and a ``special part'' decomposition (analogous to the decomposition in Theorem ~\ref{sp}). 
For more information on this axiomatic approach to closure operations, see \cite{E}.

We choose to restrict our attention to tight closure throughout the paper, since this is the most interesting operation satisfying these conditions (under certain mild assumptions on the ring) that we are aware of.
The Frobenius closure is another example of an operation satisfying these conditions (under the assumption that $(R, \m)$ is a local excellent normal ring with perfect residue field; see \cite{E2}), and therefore the result in ~\ref{gras} holds, with essentially the same proof, if one replaces ``tightly closed'' by ``Frobenius closed'', and $*$-spread by F-spread (which is equal to the number of generators of $I\q $ for $q \gg 0$).  

 Integral closure, on the other hand, while it satisfies a version of Nakayama's lemma, and it admits spread (the analytic spread), does not have a ``special part decomposition'' (see \cite{E2}). 
Therefore, Watanabe's results for integrally closed ideals cannot be recovered directly by the methods of this paper. Indeed, his methods take advantage of the geometric nature of integral closure, i.e. its connection with blowing up.

Our results on chains and families of tightly closed ideals are contained in Section 2. Section 1 deals with the connection between minimal $*$-reductions (i.e. reductions with respect to the tight closure operation) of a given ideal $I$ and the ``special part'' decomposition of $I$, giving a characterization of the minimal $*$-reductions of $I$ (Theorem~\ref{reduction}).

Throughout this paper, $(R, \m)$ will be a local excellent normal ring of characteristic $p>0$, with perfect residue field. In particular,
this implies that $R$ has a weak test element.

\section{Special tight closure and minimal $*$-reductions}
The notions of special tight closure and minimal $*$-reductions play an important role in our investigation of families of tightly closed ideals.
These notions were studied in \cite{V}, \cite{HV}, \cite{E}.
 Since the topic of this paper places us in the context of ideals $I$ that contain a given ideal $J$, we are forced to work with the slightly modified notion of minimal $*$-reduction modulo $J$.

We begin with a review of the definitions and relevant facts concerning these notions.
\begin{definition}
Let $(R, \m)$ be a local excellent normal ring of characteristic $p>0$, and let $I$ be an ideal. 
Let $R^0$ denote the set of elements in $R$ that are not in any minimal prime of $R$. 

We say that $x \in R$ is in the {\it special tight closure} of $I$ ($x \in I^{*sp}$) if there exists a $c \in R^0$ 
(equivalently, for every weak test element $c\in R$) and a fixed power of the characteristic $q_0$ such that
$cx^q \in m^{q/q_0}I\q$ for all $q=p^e\ge q_0$.

Equivalently, $x\in I^{*sp}$ if and only if there exists $q_0=p^{e_0}$ such that $x^{q_0}\in (\m I^{[q_0]})^*$.
\end{definition}

\begin{obs}(\cite{E}, Lemma 3.4)
If $K\subset I$ is such that $K^*=I^*$, then $K^{*sp}=I^{*sp}$.
\end{obs}

\begin{theorem}(\cite{HV}, Theorem 2.1)\label{sp}
Let $(R, \m)$ be an excellent normal ring of positive characteristic, with perfect residue field. Then for every ideal $I\subset R$, we have $I^*=I+I^{*sp}$.
\end{theorem}

\begin{definition}
Let $J$ be an ideal, and $f_1, \ldots, f_l \in R$. We say that $f_1, \ldots, f_l$ are 
{\it $*$-independent modulo $J$} if $f_i \notin (J, f_1, \ldots, \hat{f_i}, \ldots, f_l)^*$ for all $i=1, \ldots, l$. 

Let $J\subset I$ be ideals.
We say that $K=(J, f_1, \ldots, f_l)\subset I$ is a minimal $*$-reduction of $I$ modulo $J$ if $I\subset K^*$ and $f_1, \ldots, f_l$ are $*$-independent modulo $J$ (equivalently, $K=(J, f_1, \ldots, f_l)$ satisfies $I\subset K^*$ and is minimal with this property).
\end{definition}

\begin{obs}\label{spread}
The proof of Proposition3.3 in \cite{V} can be modified slightly to show that the property that $K=(J,f_1, \ldots, f_l)$ is generated modulo $J$ by elements which are $*$-independent modulo $J$ does not depend on the choice of a minimal system of generators for $K/J$.

The proof of Theorem 5.1 in \cite{E} can be modified slightly to show that if $K, K'$ are minimal $*$-reductions for $I$ modulo $J$, then $K/J$ and $K'/J$ have the same minimal number of generators, which will be called the $*$-{\it spread of} $I$ {\it modulo} $J$, and will be denoted $l$, or $l^*_J(I)$.
\end{obs}

We also recall two important properties of tight closure that will be required in the proofs of our results:

\begin{prop}(\cite{A}, Proposition 2.4)\label{colon}
Let $(R, \m)$ be an excellent analytically irreducible local ring of characteristic $p>0$, let $I$ be an ideal, and let $f \in R$. 

Assume $f \notin I^*$; then there exists $q_0=p^{e_0}$ such that, for all $q\ge q_0$, we have $I\q :f^q \subset \m^{q/q_0}$.
\end{prop}

Note that our standard hypothesis ($(R, \m)$ is a local excellent normal ring) implies that $R$ is analytically irreducible, since the completion of an excellent normal ring is again normal, hence a domain.

The following is known as the ``Nakayama lemma for tight closure'':
\begin{prop}(\cite{E}, Proposition 2.1)\label{nak}
Let $(R, \m)$ be a Noetherian local ring of characteristic $p>0$, possessing a $q_0$-weak test element $c$. Let $I, J$ be ideals of $R$ such that $J\subset I\subset (J+\m I)^*$. Then $I\subset J^*$.
\end{prop}

The following observation follows immediately from Proposition ~\ref{nak}:
\begin{corollar}\label{min}
If $K=(J, f_1, \ldots, f_l)$ is a minimal $*$-reduction of $I$ modulo $J$, then the images of $f_1, \ldots, f_l$ must be part of a minimal system of generators for $I/J$.
\end{corollar}
\begin{proof}
Assume that $f_i \in (J, \m I, f_1, \ldots, f_{i-1}, f_{i+1}, \ldots, f_l)$ for some $i \in \{1, \ldots, l\}$ and seek a contradiction.
Then
$$I^*\subset(J, f_1, \ldots, f_l)^*\subset (J, \m I, f_1, \ldots, f_{i-1}, f_{i+1}, \ldots, f_l)^*,$$
and Proposition~\ref{nak} shows that
$I^*\subset (J, f_1, \ldots, f_{i-1}, f_{i+1}, \ldots, f_l)^*$, contradicting the $*$-independence, modulo $J$, of $f_1, \ldots, f_l$.
\end{proof}

 We will also use the following observation, which follows immediately from the definition of tight closure:
\begin{obs}\label{obvious}

a). If $f^{q_0}\in (I^{[q_0]})^*$ for some $q_0=p^{e_0}$, then $f \in I^*$.

b). $(I_1+I_2^*)^*=(I_1+I_2)^*$ for any ideals $I_1$, $I_2$.
\end{obs}

The main result of this section is a characterization of the ideals $K\subset I$ which are minimal $*$-reductions for $I$ modulo $J$ (Theorem ~\ref{reduction}). We also consider the question of which ideals $K=(J, f_1,\ldots, f_i)$ can be extended to minimal $*$-reductions (Proposition ~\ref{extend}). These results will play an important role in our characterization of the family of adjacent tightly closed ideals in the next section.

We begin with a preliminary result illustrating the connection between special tight closure and $*$-reductions. This is a generalization of Theorem 2.1 in \cite{HV}:

\begin{prop}\label{sum}
Let $J\subset K\subset I$ be ideals of $R$ with $K$ a minimal $*$-reduction of $I$ modulo $J$.
Then we have a direct sum decomposition
$$
\frac{I^*}{\m I+J}\cong \frac{\m I+K}{\m I + J} \oplus \frac{I^{*sp}+J}{\m I +J}.
$$ 
\end{prop}

\begin{obs}\label{exclude}
By Theorem~\ref{sp} we know that $I^*=K^*=K+K^{*sp}=K+I^{*sp}$.
The new content of Proposition~\ref{sum} is the fact that
$
I^{*sp}\cap K\subset (\m I+J).
$
Equivalently, 
if $g \in K\, \backslash \, (\m I +J)$, then $g \notin I^{*sp}$.
\end{obs}
\begin{proof}
It is enough to prove the equivalent formulation given in the last paragraph of Observation~\ref{exclude}.

We may extend $g$ to a minimal system of generators $g=g_1, \ldots, g_l$ for $K/J$; note that $g_1, \ldots, g_l$ may be chosen to be $*$-independent modulo $J$.

 Assume that $g \in I^{*sp}$ and seek 
a contradiction; thus there exists $q_0=p^{e_0}$ such that $g^{q_0}\in (\m I^{[q_0]})^*$. It follows that
$$
(I^{[q_0]})^*\subset (J^{[q_0]}, g_1^{q_0}, \ldots, g_l^{q_0})^*\subset
(J^{[q_0]}, g_2^{q_0}, \ldots, g_l^{q_0}, (\m I^{[q_0]})^*)^*.
$$
Using ~\ref{obvious}(b), Proposition~\ref{nak} implies that $(I^{[q_0]})^*\subset (J^{[q_0]}, g_2^{q_0}, \ldots, g_l^{q_0})^*$, and thus, by 
~\ref{obvious}(a), we get $I^*\subset (J, g_2, \ldots, g_l)^*$, contradicting the $*$-independence of $g_1, \ldots, g_l$ modulo $J$.
\end{proof}

The next result deals with the following question: under what circumstances can one obtain a new minimal $*$-reduction from an existing one, by replacing one of the generators and keeping the others?

\begin{prop}\label{switch}
Let $J\subset I$ be ideals, and let $f \in R$.

a). Assume that $R$ is analytically irreducible. Let $K=(J, f_1, \ldots, f_l)$ be a minimal $*$-reduction for $I$ modulo $J$. 

If $f\in I\, \backslash\, I^{*sp}$, then there exists an $i \in \{1, \ldots, n\}$ such that $K'=(J,  f_1, \ldots, f_{i-1}, f, f_{i+1}, \ldots, f_l)$ is again a minimal $*$-reduction for $I$ modulo $J$.

b). If $f_1, \ldots, f_l \in I$ are $*$-independent modulo $J$, then they are also $*$-independent modulo $J+I^{*sp}$.

 In particular, if $K=(J, f_1, \ldots, f_l)$ is a minimal $*$-reduction for $I$ modulo $J$ and $f\in I^{*sp}$, 
then $K'=(J,  f_1, \ldots, f_{i-1},f, f_{i+1}, \ldots, f_l)$ is not a minimal $*$-reduction for $I$ modulo $J$, 
for any choice of $i$.
\end{prop}
\begin{proof}
a). According to Theorem~\ref{sp}, there exists $f'=f+\alpha_1 f_1 +\ldots + \alpha _l f_l \in I^{*sp}$.
Since $\m I \subset I^{*sp}$ and $f \notin I^{*sp}$, at least one $\alpha _i$ is not in $\m$.
Choose $c \in R^{0}$ and $q_0=p^{e_0}$ such that for all $q=p^e>q_0$ we have
$cf'^q =b_1f_1^q +\ldots +b_lf_l^q$ (mod $J\q$), with $b_1, \ldots, b_l \in m^{q/q_0}$.
It follows that 
$$(c\alpha _i^q-b_i)f_i^q \in (J\q, f^q, f_1^q, \ldots, f_{i-1}^q, f_{i+1}^q, \ldots, f_l^q).$$ 
Theorem ~\ref{colon} implies
$f_i \in (J, f, f_1, \ldots, f_{i-1}, f_{i+1}, \ldots, f_l)^*$, and therefore 
$I^*\subset (J, f, f_1, \ldots, f_{i-1}, f_{i+1}, \ldots, f_l)^*$. This is enough to imply that
$$(J, f, f_1, \ldots, f_{i-1}, f_{i+1}, \ldots, f_l)$$ is a minimal $*$-reduction, for otherwise one of $f, f_1, \ldots, f_{i-1}, f_{i+1}, \ldots, f_l$ would be redundant, and we would be able to find a minimal $*$-reduction with fewer generators.

b). For each $i \in \{1, \ldots, l\}$, we need to show that 
$$f_i \notin (J, I^{*sp}, f_1, \ldots, f_{i-1}, f_{i+1}, \ldots, f_l)^*.$$ Otherwise, one could extract a minimal $*$-reduction $K''$ modulo $J$ generated by some of $f_1, \ldots, f_{i-1}, f_{i+1}, \ldots, f_l,$ and some of the elements in $I^{*sp}$. Since elements among
$f_1, \ldots, f_{i-1}, f_{i+1}, \ldots, f_l$ alone cannot generate a  $*$-reduction, it follows that we 
must have elements in $I^{*sp}$ among the minimal generators of $K''/J$. But $I^{*sp}\cap K''\subset \m I+J$ by Observation~\ref{exclude}, and this is a contradiction.
\end{proof}

\begin{theorem}\label{reduction}
Let $J\subset K\subset I$ be ideals of $R$ such that $K$ is generated by $l$ elements which are part of a minimal system of generators for $I/J$, where $l=l^*_J(I)$. Assume that $I$ is tightly closed.

The following are equivalent:

a. $K$ is a minimal $*$-reduction for $I$ modulo $J$;

b. $I^{*sp}\cap K\subset \m I+J$;

c. $I=I^{*sp}+K$.
\end{theorem}
\begin{proof}
a. $\Rightarrow $ b. This is contained in Observation ~\ref{exclude}.

b. $\Leftrightarrow $ c. 
Choose $K'=(J, f_1, \ldots, f_l)$ a minimal $*$-reduction for $I$ modulo $J$. 
 Proposition~\ref{sum} shows that 
$$
\mathrm{dim}\left(\frac{I}{\m I+J}\right)=\mathrm{dim}\left(\frac{I^{*sp}+J}{\m I+J}\right)+\mathrm{dim}\left(\frac{\m I+K'}{\m I+J}\right)= \mathrm{dim}\left(\frac{I^{*sp}+J}{\m I+J}\right)+l
$$
(these are dimensions as vector spaces over $k=R/\m$).

 Since dim$((\m I +K)/(\m I+J))=l$, the equivalence follows.

c. $\Rightarrow $ a.
Choose $q_0\gg 0$ so that we have
$$
I^{[q_0]}=(I^*)^{[q_0]}=K^{[q_0]}+(I^{*sp})^{[q_0]}\subset ( K^{[q_0]}+(\m I^{[q_0]})^*)^*=(K^{[q_0]}+\m I^{[q_0]})^*
$$
By Proposition ~\ref{nak}, this implies $I^{[q_0]}\subset (K^{[q_0]})^*$, and therefore $I\subset K^*$. This is sufficient to show that $K$ is a minimal $*$-reduction, since otherwise we could find a minimal $*$-reduction with fewer generators.
\end{proof}

We are also interested in characterizing the ideals which can be extended to minimal $*$-reductions by adding extra generators:

\begin{corollar}\label{extend}
Let $J\subset  I$ with $I$ tightly closed, and let $i\le l=l^*_J(I)$. Let $g_1, \ldots, g_i$ be part of a minimal system of generators for $I$ modulo $J$
and let $K=(J, g_1,\ldots, g_i)$. 

The following are equivalent:

a. $K$ can be extended to a minimal $*$-reduction of $I$ modulo $J$, i.e. 
there exist $g_{i+1}, \ldots g_l \in I$ such that $(J, g_1,\ldots, g_l)$ is a minimal $*$-reduction for $I$ modulo $J$.

b. $I^{*sp}\cap K\subset \m I+J$
\end{corollar}
\begin{proof}
a. $\Rightarrow$ b. This follows immediately from Theorem ~\ref{reduction}.

b. $\Rightarrow $ a. 
We need to prove that 
 we can find $g_{i+1}, \ldots g_l$ in $I/(\m I+J)$
such that $g_1, \ldots, g_l$ are linearly independent in $I/(\m I+J)$, and 
$K'=(J, g_1, \ldots, g_l)$ satisfies the condition in b. of Theorem ~\ref{reduction}.
This can be achieved by extending the images of $g_1, \ldots, g_i $ in $I/(\m I+J)$ to a basis for a subspace 
containing $(\m I+ K)/(\m I+J)$ and complementary to $(I^{*sp}+J)/(\m I+J)$.
\end{proof}

\section{Chains and families of tightly closed ideals}

\begin{definition}
Let $J\subset I$ be tightly closed ideals. Let ${\mathcal F}(J, I)$ be the set of all tightly closed ideals $I'$ such that $J\subset I'\subset I$, and $\lambda(I/I')=1$. 
\end{definition}

\begin{theorem}\label{family}
${\mathcal F}(J, I)$ is the set of all ideals $I'$  of the form
$$
I'=(J, f_1, \ldots, f_{l-1})+ I^{*sp}
$$
where 
$f_1, \ldots, f_{l-1}\in I$ are such that $(J, f_1, \ldots, f_{l-1})$ can be extended to a minimal $*$-reduction $K=(J, f_1, \ldots, f_l)$ of $I$ modulo $J$, for some choice of $f_l \in I$.
\end{theorem}
\begin{proof}
First we prove that any ideal $I'$ of the given form is in ${\mathcal F}(J, I)$. Since $I=K+I^{*sp}$ and $\m I\subset I^{*sp}$, it follows that $\lambda(I/I')=1$, and in fact $I/I'$ is spanned by the image of $f_l$. In order to see that $I'$ is tightly closed, note that we have $f_l \notin I'^*$ by Proposition ~\ref{switch} part b). On the other hand, $I'^*\subset I^*=I$, and every element in $I\, \backslash I'$ is congruent modulo $I'$ to a unit multiple of $f_l$. This shows that $I'$ is tightly closed.

Conversely, we need to show that every ideal $I'\in {\mathcal F}(J, I)$ has the given form. To this end, it suffices to show that there exist $f_1, \ldots, f_{l-1}$ such that $(J, f_1, \ldots, f_l)$ is a minimal $*$-reduction for some choice of $f_l$, and 
$(J, f_1, \ldots, f_{l-1})+ I^{*sp}\subset I'$.
Once we have this inclusion, the equality follows, since
$$
\mathrm{dim}\left(\frac{I'}{\m I+J}\right)=\mathrm{dim}\left(\frac{I}{\m I+J}\right)-1=l-1+\mathrm{dim}\left(\frac{I^{*sp}+J}{\m I+J}\right),
$$
and Corollary ~\ref{extend} shows that this is equal to the dimension of  
$$\frac{(J, f_1, \ldots, f_{l-1})+ I^{*sp}}{\m I+J}.
$$

Choose $K=(J, f_1, \ldots, f_l)$ an arbitrary minimal $*$-reduction for $I$ modulo $J$. Since $I'$ is tightly closed, 

$K\not\subset  I'$. Since $\lambda(I/I')=1$, it follows that $I'+K=I$, and therefore $\lambda(K/I'\cap K)=\lambda((I'+K)/I')=1$. 

This implies that we can choose generators for $K/J$ such that $(J, f_1, \ldots, f_{l-1}, \m f_l)\subset I'$, $f_l \notin I'$. 

Therefore, the image of $f_l$ generates $I/I'$. It remains to be shown that $I^{*sp}\subset I'$. Let $g\in I^{*sp}$. Assume
that $g \notin I'$ and seek a contradiction. Then $g \equiv uf_l$ (mod $I'$), with $u \in R$ a unit (since $\m I\subset I'$).
Choose $q_0$ such that $g^{q_0}\in (\m I^{[q_0]})^*$. We have
$$
I^{[q_0]}\subset (J^{[q_0]}, f_1^{q_0}, \ldots, f_l^{q_0})^*= (J^{[q_0]}, f_1^{q_0}, \ldots, f_{l-1}^{q_0}, g^{q_0}, I'^{[q_0]})^*\subset
(I'^{[q_0]},  \m  I^{[q_0]})^*.
$$
The last inclusion follows since $(J, f_1, \ldots, f_{l-1})\subset I'$.
Proposition ~\ref{nak} implies $I^{[q_0]}\subset (I'^{[q_0]})^*$, and thus $I\subset I'^*$, which is a contradiction.
\end{proof}

\begin{corollar}\label{chain}
The family ${\mathcal F}(J, I)$ is non-empty. In particular, if $I, J$ are $\m$-primary, there exists a sequence $J=I_0\subset I_1 \ldots \subset I_n=I$ consisting of tightly closed ideals, such that $\lambda(I_{i+1}/I_i)=1$ for all $i$.
\end{corollar}

\begin{corollar}\label{gras}
The ideals in ${\mathcal F}(J, I)$ are in one-to-one correspondence with points on the Grassmanian variety of $l-1$ dimensional subspaces $(J, f_1, \ldots, f_{l-1})+ I^{*sp})/(J+I^{*sp})$ in the $l$-dimensional vector space $V=I/(J+I^{*sp})$, where $l$ denotes the $*$-spread of $I$ modulo $J$.

In particular, the ideals in ${\mathcal F}(J, I)$ are in one-to-one correspondence with points of an $l-1$ dimensional projective space ${\bf P}(V^*)$.
\end{corollar}
\begin{proof}
Note that $V=I/(J+I^{*sp})$ is a vector space because $\m I\subset I^{*sp}$, and we can write it as $(K+I^{*sp})/(J+I^{*sp})$, where $K$ is a minimal $*$-reduction for $I$ modulo $J$. According to Theorem~\ref{reduction}, we have
a direct sum decomposition
$$
\frac{K+I^{*sp}}{\m I+J}=\frac{I^{*sp}+J}{\m I+J}\oplus \frac{K+\m I}{\m I+J};
$$
therefore we have
$$
\mathrm{dim}\left(\frac{K+I^{*sp}}{\m I+J}\right) =\mathrm{dim}\left(\frac{I^{*sp}+J}{\m I+J}\right)+l,
$$
and this implies 
$$
\mathrm{dim}(V)=\mathrm{dim}\left(\frac{K+I^{*sp}}{\m I+J}\right)-\mathrm{dim}\left(\frac{I^{*sp}+J}{\m I+J}\right)=l.
$$

According to Theorem~\ref{family}, the set ${\mathcal F}(J, I)$ is in one-to-one correspondence with subspaces
 of the form $$\frac{(J, f_1, \ldots, f_{l-1})+I^{*sp}}{J+I^{*sp}}$$ of $V$, where $K'=(J, f_1, \ldots, f_{l-1})$ satisfies the condition in Corollary ~\ref{extend}. 

This implies that
$$
\mathrm{dim}\left(\frac{K'+I^{*sp}}{\m I+J}\right) =\mathrm{dim}\left(\frac{I^{*sp}+J}{\m I+J}\right)+\mathrm{dim}\left(\frac{K'+\m I}{\m I+J}\right) = 
$$
$$
\mathrm{dim}\left(\frac{I^{*sp}+J}{\m I+J}\right)+ l-1,$$
from which it follows that
$$
\mathrm{dim}\left( \frac{K'+I^{*sp}}{I^{*sp}+J}\right)=l-1.
$$

Conversely, every $l-1$ dimensional subspace $W$ of $V$ is of this form
$(K'+I^{*sp})/(J+I^{*sp})$ with $K'=(J, f_1, \ldots, f_{l-1})$ satisfying the equalities above, which are equivalent to the condition in Corollary ~\ref{extend}.
\end{proof}

\bigskip

{\sf Department of Mathematics, University of South Carolina, Columbia, SC}

\bigskip

{\it E-mail address:} vraciu@math.sc.edu

\end{document}